\documentclass[12pt,letter]{article}

\usepackage{mystyle}
\usepackage{mythmstyle_report}

\begin{document}

\title{Automorphism group of the modified bubble-sort graph}
\author{Ashwin~Ganesan%
  \thanks{Department of Electronics and Telecommunication Engineering, Vidyalankar Institute of Technology, Wadala, Mumbai, India. Correspondence address: 
\texttt{ashwin.ganesan@gmail.com}}
}

\date{}

\maketitle

\begin{abstract} 
The modified bubble-sort graph of dimension $n$ is the Cayley graph of $S_n$ generated by $n$ cyclically adjacent transpositions.   In the present paper, it is shown that the automorphism group of the modified bubble sort graph of dimension $n$ is $S_n \times D_{2n}$, for all $n \ge 5$.  Thus, a complete structural description of the automorphism group of the modified bubble-sort graph is obtained.  A similar direct product decomposition is seen to hold for arbitrary normal Cayley graphs generated by transposition sets.
\end{abstract}
\bigskip

\section{Introduction}

Let $X=(V,E)$ be a simple undirected graph.  The (full) automorphism group of $X$, denoted by $\Aut(X)$, is the set of permutations of the vertex set that preserves adjacency, i.e., $\Aut(X):=\{g \in \Sym(V): E^g=E\}$.  Let $H$ be a group with identity element $e$, and let $S$ be a subset of $H$.  The Cayley graph of $H$ with respect to $S$, denoted by $\Cay(H,S)$, is the graph with vertex set $H$ and arc set $\{ (h,sh): h \in H, s \in S\}$.  When $S$ satisfies the condition  $1 \notin S = S^{-1}$, the Cayley graph $\Cay(H,S)$ has no self-loops and can be considered to be undirected.

A Cayley graph $\Cay(H,S)$ is vertex-transitive since the right regular representation $R(H)$ acts as a group of automorphisms of the Cayley graph. The set of automorphisms of $H$ that fixes $S$ setwise is a subgroup of the stabilizer $\Aut(\Cay(H,S))_e$ (cf. \cite{Biggs:1993}, \cite{Godsil:Royle:2001}).  A Cayley graph $X:=\Cay(H,S)$ is said to be \emph{normal} if $R(H)$ is a normal subgroup of $\Aut(X)$, or equivalently, if $\Aut(X) = R(H) \rtimes \Aut(H,S)$ (cf. \cite{Xu:1998}).     

Let $S$ be a set of transpositions generating the symmetric group $S_n$.  The transposition graph of $S$, denoted by $T(S)$, is defined to be the graph with vertex set $\{1,\ldots,n\}$, and with two vertices $i$ and $j$ being adjacent in $T(S)$ whenever $(i,j) \in S$.  A set $S$ of transpositions generates $S_n$ iff the transposition graph of $S$ is connected.   When the transposition graph of $S$ is the $n$-cycle graph, then the Cayley graph $\Cay(S_n,S)$ is called the modified bubble-sort graph of dimension $n$.  Thus, the modified bubble-sort graph of dimension $n$ is the Cayley graph of $S_n$ with respect to the set of generators $\{(1,2),(2,3),\ldots,(n-1,n),(n,1)\}$.  The modified bubble-sort graph has been investigated for consideration as the topology of interconnection networks (cf. \cite{Lakshmivarahan:etal:1993}). Many authors have investigated the automorphism group of graphs that arise as the topology of interconnection networks; for example, see \cite{Deng:Zhang:2011}, \cite{Deng:Zhang:2012}, \cite{Ganesan:JACO}, \cite{Zhang:Huang:2005}, \cite{Zhou:2011}. 

Godsil and Royle \cite{Godsil:Royle:2001} showed that if the transposition graph of $S$ is an asymmetric tree, then the automorphism group of the Cayley graph $\Cay(S_n,S)$ is isomorphic to $S_n$.  Feng \cite{Feng:2006} showed that $\Aut(S_n,S)$ is isomorphic to $\Aut(T(S))$ and that if the transposition graph of $S$ is an arbitrary tree, then the automorphism group of $\Cay(S_n,S)$ is the semidirect product $R(S_n) \rtimes \Aut(S_n,S)$.  Ganesan \cite{Ganesan:DM:2013} showed that if the girth of the transposition graph of $S$ is at least 5, then the automorphism group of the Cayley graph $\Cay(S_n,S)$ is the semidirect product $R(S_n) \rtimes \Aut(S_n,S)$.  The results in the present paper imply that all these automorphism groups in the literature can be factored as a direct product.

In Zhang and Huang \cite{Zhang:Huang:2005}, it was shown the automorphism group of the modified bubble-sort graph of dimension $n$ is the group product $S_n D_{2n}$ (groups products are also referred to as Zappa-Szep products).  This result was strengthened in Feng \cite{Feng:2006}, where it was proved that the automorphism group of the modified bubble-sort graph of dimension $n$ is the semidirect product $R(S_n) \rtimes D_{2n}$ (cf. \cite[p. 72]{Feng:2006} for an explicit statement of this conclusion).  

In the present paper, we obtain a complete structural description of the automorphism group of the modified bubble-sort graph of dimension $n$:

\begin{Theorem} \label{theorem:aut:mbs:dp}
 The automorphism group of the modified bubble-sort graph of dimension $n$ is $S_n \times D_{2n}$, for all $n \ge 5$.
\end{Theorem}

We shall prove the following more general result:

\begin{Theorem} \label{thm:main:dp}
 Let $S$ be a set of transpositions generating $S_n (n \ge 3)$ such that the Cayley graph $\Cay(S_n,S)$ is normal.  Then, the automorphism group of the Cayley graph $\Cay(S_n,S)$ is the direct product $S_n \times \Aut(T(S))$, where $T(S)$ denotes the transposition graph of $S$.
\end{Theorem}

In the special case where $T(S)$ is the $n$-cycle graph, $\Aut(T(S))$ is isomorphic to the dihedral group $D_{2n}$ of order $2n$.  Hence, Theorem~\ref{theorem:aut:mbs:dp} is a special case of Theorem~\ref{thm:main:dp}. Also, Ganesan \cite{Ganesan:DM:2013} showed that the modified bubble-sort graphs of dimension less than 5 are non-normal; hence, the assumption $n \ge 5$ in Theorem~\ref{theorem:aut:mbs:dp} is necessary.

\bigskip \noindent \emph{Remark 1}. Given a set $S$ of transpositions generating $S_n$, let $G:=\Aut(\Cay(S_n,S))$.  In the instances where $G = R(S_n) \rtimes G_e$, the factor $G_e \cong \Aut(T(S))$ is in general not a normal subgroup of $G$, and so the semidirect product cannot be written immediately as a direct product.  For example, for the modified bubble-sort graph of dimension $n$, $G \cong R(S_n) \rtimes G_e \cong S_n \rtimes D_{2n}$, where $G_e$ is not normal in $G$.  In the present paper, it is shown that $R(S_n)$ has another complement in $G$ which is a normal subgroup of $G$.   In the proof below, we show that the image of $\Aut(T(S))$ under the left regular action of $S_n$ on itself is a normal complement of $R(S_n)$ in $G$.  Thus, the direct factor $\Aut(T(S))$ that arises in $G \cong R(S_n) \times \Aut(T(S))$ is not $G_e$ but is obtained in a different manner.


\section{Proof of Theorem \ref{thm:main:dp}}

Let $S$ be a set of transpositions generating $S_n$.  We first establish that the Cayley graph $\Cay(S_n,S)$ has a particular subgroup of automorphisms.   In this section, let $\lambda$ denote the left regular action of $S_n$ on itself, defined by  $\lambda: S_n \rightarrow \Sym(S_n), a \mapsto \lambda_a$, where $\lambda_a: x \mapsto a^{-1}x$.  

\begin{Proposition} \label{prop:leftaction:auts}
Let $T(S)$ denote the transposition graph of $S$. Then, $\{ \lambda_a: a \in \Aut(T(S)) \}$ is a set of automorphisms of the Cayley graph $X:=\Cay(S_n,S)$.
\end{Proposition}

\noindent \emph{Proof}:
Let $a \in \Aut(T(S))$.  We show that $\{h,g\} \in E(X)$ if and only if $\{h,g\}^{\lambda(a)} \in E(X)$.  Suppose $\{h,g\} \in E(X)$.  Then $g=sh$ for some transposition $s=(i,j) \in S$.  We have that $\{h,g\}^{\lambda(a)}  = \{h,sh\}^{\lambda(a)} = \{h^{\lambda(a)},(sh)^{\lambda(a)}\} =  \{a^{-1} h, a^{-1}sh\} = \{a^{-1} h, (a^{-1} s a) a^{-1} h\}$. Now $a^{-1} s a = a^{-1} (i,j) a =  (i^a,j^a) \in S$ since $a$ is an automorphism of the graph $T(S)$ that has edge set $S$. Thus, $\{h,sh\}^{\lambda(a)} \in E(X)$. Conversely,  suppose $\{h,g\}^{\lambda(a)} \in E(X)$.  Then $a^{-1}h = s a^{-1}g$ for some $s \in S$. Hence $h = (a s a^{-1}) g$.  We have that $a s a^{-1} = a(i,j)a^{-1}=(i,j)^{a^{-1}} \in S$ because $a$ is an automorphism of $T(S)$.  Hence, $h$ is adjacent to $g$.  Thus, $\lambda(\Aut(T(S)))$ is a subgroup of $\Aut(X)$.  
\qed

\begin{Theorem} 
 Let $S$ be a set of transpositions generating $S_n (n \ge 3)$ such that the Cayley graph $\Cay(S_n,S)$ is normal.  Then, the automorphism group of the Cayley graph $\Cay(S_n,S)$ is $S_n \times \Aut(T(S))$, where $T(S)$ denotes the transposition graph of $S$.
\end{Theorem}

\noindent \emph{Proof}: Let $X$ denote the Cayley graph $\Cay(S_n,S)$.  Since $X$ is a normal Cayley graph, its automorphpism group $\Aut(X)$ is equal to $R(S_n) \rtimes \Aut(S_n,S)$ (cf. \cite{Xu:1998}).   Let $R(a)$ denote the permutation of $S_n$ induced by action by right multiplication by $a$, so that $R(S_n):=\{R(a): a \in S_n\}$ is the right regular representation of $S_n$.    The intersection of the left and right regular representations of a group is the image of the center of the group under either action.  The center of $S_n$ is trivial, whence $R(S_n) \cap \lambda(S_n) =1$.  In particular, $\lambda(\Aut(T(S)))$ and $R(S_n)$ have a trivial intersection.  By Feng \cite{Feng:2006}, $\Aut(S_n,S) \cong \Aut(T(S))$, and it follows from cardinality arguments that $R(S_n) \lambda(\Aut(T(S)))$ exhausts all the elements of $\Aut(X)$. Thus, $R(S_n)$ and $\lambda(\Aut(T(S)))$ are complements of each other in $\Aut(X)$ and every element in $\Aut(X)$ can be expressed uniquely in the form $R(a) \lambda(b)$ for some $a \in S_n$ and $b \in \Aut(T(S))$.  This proves that $\Aut(X) = R(S_n) \rtimes \lambda(\Aut(T(S)))$. 

It remains to prove that $\lambda(\Aut(T(S)))$ is a normal subgroup of $\Aut(X)$.   Suppose $g \in \Aut(X)$ and $c \in \Aut(T(S))$. We show that $g^{-1} \lambda(c) g \in \lambda(\Aut(T(S)))$. We have that $g = R(a) \lambda(b)$ for some $a \in S_n, b \in \Aut(T(S))$.  Hence, $g^{-1} \lambda(c) g = (R(a) \lambda(b))^{-1} \lambda(c) (R(a) \lambda(b))$, which maps $x \in S_n$ to $b^{-1} c^{-1} bxa^{-1}a = b^{-1} c^{-1} bx$.  Since $b,c \in \Aut(T(S))$,  $d^{-1}:=b^{-1}c^{-1}b \in \Aut(T(S))$. Thus, $g^{-1} \lambda(c) g = \lambda(d) \in \lambda(\Aut(T(S)))$. Hence, $\lambda(\Aut(T(S)))$ is a normal subgroup of $\Aut(X)$ and $\Aut(X) = R(S_n) \times \lambda(\Aut(T(S)))$.  Since $\lambda(\Aut(T(S))) \cong \Aut(T(S))$, the assertion follows.
\qed

\bigskip \noindent \emph{Remark 2.} We recall a particular result from group theory, which can be used to deduce that the semidirect products in the literature can be strengthened to direct products. 
Let $A$ be a subgroup of a group $H$ and suppose $H$ has a trivial center. 
 Let $A$ act on $H$ by conjugation.  Let $\lambda(A)$ denote the image of the left action of $A$ on $H$.
Then the groups $R(H) \rtimes \Inn(A)$ and $R(H) \times \lambda(A)$ are isomorphic, where both groups are internal group products and subgroups of $\Sym(H)$.  It follows from this group-theoretic result that the automorphism group of the Cayley graphs mentioned above can be factored as direct products. 
However, to the best of our knowledge, this group-theoretic result has not been used so far to deduce results in the context of automorphism groups of Cayley graphs generated by transposition sets - the expressions given in the previous literature for the automorphism group of Cayley graphs mentioned above have been only semidirect product factorizations (cf. \cite[p.72]{Feng:2006}, \cite{Ganesan:DM:2013}, \cite{Xu:1998}). In the present paper, in addition to obtaining a complete structural description of the automorphism group of the modified bubble-sort graph and of a family of normal Cayley graphs, the proof method also includes Proposition~\ref{prop:leftaction:auts}, which establishes that these graphs possess certain automorphisms.

\bibliographystyle{plain}
\bibliography{refsaut}


\end{document}